\documentclass[11pt]{article}
\usepackage{latexsym,amsmath,amssymb,amsfonts,epsfig,graphicx}
\usepackage{latexsym,amsmath,amssymb,amsfonts,epsfig,graphicx,cite,psfrag}
\usepackage{eepic,color,colordvi,amscd,amsthm}

\newtheorem{theo}{Theorem}[section]
\newtheorem{lem}{Lemma}[section]

\newtheorem{conj}{Conjecture}

\def\qed{\hfill \rule{4pt}{7pt}}
\def\pf{\noindent {\it Proof.} }

\begin{document}

\title{Factor-Critical Property in $3$-Dominating-Critical Graphs
\thanks{This work is supported by RFDP of Higher Education of China and Natural Sciences and
Engineering Research Council of Canada.}}

\author{Tao Wang$^1$  and Qinglin Yu$^1$$^2$$\thanks{  {\it Corresponding author:
yu@tru.ca}}$
\\ {\small  $^1$Center for Combinatorics, LPMC}
\\ {\small Nankai University, Tianjin, China}
\\ {\small  $^2$Department of Mathematics and Statistics}
\\ {\small Thompson Rivers University, Kamloops, BC, Canada}
}
\maketitle

\begin{abstract}
    A vertex subset $S$ of a graph $G$
is a {\it dominating set} if every vertex of $G$ either belongs to
$S$ or is adjacent to a vertex of $S$. The cardinality of a
smallest dominating set is called the {\it dominating number} of
$G$ and is denoted by $\gamma(G)$. A graph $G$ is said to be
$\gamma$-{\it vertex-critical} if $\gamma(G-v)< \gamma(G)$, for every
vertex $v$ in $G$.

    Let $G$ be a $2$-connected $K_{1,5}$-free $3$-vertex-critical graph.
For any vertex $v \in V(G)$, we show that $G-v$ has a perfect matching
(except two graphs), which is a
conjecture posed by Ananchuen and Plummer \cite{PLU2}.\begin{flushleft}
\end{flushleft}
\end{abstract}

\section{Introduction}
Let $G$ be a finite simple graph with vertex set $V(G)$ and edge
set $E(G)$. A vertex subset $S$ of $G$ is a {\it dominating
set} of $G$ if every vertex of $G$ either belongs to $S$ or is
adjacent to a vertex of $S$. The minimum size of such a set is
called the \textit{dominating number} of
$G$ and is denoted by $\gamma(G)$.
A graph G is {\it vertex domination-critical}, or
\textit{$\gamma$-vertex-critical},  if for any vertex $v$ of $G$,
$\gamma(G-v)< \gamma(G)$. We use $G[S]$ to denote the subgraph
induced by $S$ for some $S\subseteq V(G)$. The minimum degree of $G$
is denoted by $\delta(G)$. A graph is called
\textit{$K_{1,k}$-free} if it has no induced subgraph isomorphic to
the complete bipartite graph $K_{1,k}$.

  A matching is \textit{perfect} if it is
incident with every vertex of $G$. If $G-v$ has a perfect matching,
for every choice of $v\in V(G)$, $G$ is said to be
\textit{factor-critical}.  The concept of factor-critical graphs was first introduced by
Gallai in 1963 and it plays an important role in the study of matching theory.
To be contrary to its apparent strong property, such graphs form a relatively rich
family for study. It is the essential ``building block" for
the so-called Gallai-Edmonds structure for matchings.

The subject of $\gamma$-vertex-critical graphs was studied first by
Brigham, Chinn and Dutton \cite {BCD1,BCD2} and continued by Fulman {\it et al.}
\cite {FUL1, FUL2}. Clearly, the only $1$-vertex-critical graph is
$K_1$ (a single vertex). Brigham, Chinn and Dutton \cite {BCD1}
pointed out that the $2$-vertex-critical graphs are precisely the family of
graphs obtained from the complete graphs $K_{2n}$ with a perfect
matching removed. For $\gamma >2$, however, much remains unknown
about the structure of $\gamma$-vertex-critical graphs. Recently,
Ananchuen and Plummer \cite {PLU1, PLU2} began to study matchings in
$3$-vertex-critical graphs. They showed that a $K_{1,5}$-free $3$-vertex-critical
graph of even order has a perfect matching (see \cite {PLU1}) and
a $K_{1,4}$-free $3$-vertex-critical
graph of odd order is factor-critical (see \cite {PLU2}).  Furthermore,
they posed the following conjecture.

\begin{conj}   \label{conj}
If $G$ is a $K_{1,5}$-free $3$-vertex-critical $2$-connected graph
of odd order with $\delta(G)\geq 3$, then $G$ is factor-critical.
\end{conj}

In this paper, we show that the conjecture holds for almost all
graphs and there are only two counterexamples. \\

    If $v\in V(G)$, we denote by $D_v$, a minimum dominating
set of $G-v$. The following facts about $D_v$ follow immediately from the definition
of $3$-vertex-criticality and we shall use it frequently in the proof of the main theorem.

\noindent\textbf{Facts:} If $G$ is $3$-vertex-critical, then the
followings hold:
\begin{itemize}
\vspace{-5pt}
\item [(1)] For every vertex $v$ of $G$, $|D_v|=2$.
\vspace{-5pt}
\item [(2)] If $D_v=\{x,y\}$, then $x$ and $y$ are not adjacent to $v$.
\vspace{-5pt}
\item [(3)] For every pair of distinct vertices $v$ and $w$, $D_v\neq D_w$.
\end{itemize}

  The readers are referred to \cite{LOV} for other terminology not specified in this paper.

\section{Main Result}

By Tutte's well-known $1$-Factor Theorem, if a graph $G$ has no
perfect matching, then there exits a set $S\subseteq V(G)$ such that
the number of components in $G-S$ having odd order is greater than
the order of $S$. If $S\subseteq V(G)$, we shall denote by
$\omega(G-S)$, the number of components of $G-S$ and by $c_o(G-S)$,
the number of odd components of $G-S$. A criterion similar
to $1$-Factor Theorem for factor-critical graphs is as
follows.

\begin{lem}    \label{lem1} (see \cite{LOV}) A graph $G$ is factor-critical
if and only if $c_o(G-S)\leq |S|-1$,
for every nonempty set $S\subseteq V(G)$.
\end{lem}

\begin{lem}    \label{lem2}
Let $G$ be $3$-vertex-critical and $S$ be a cutset in $G$ with
$|S|\geq 4$. If $D_u\subseteq S$ for each vertex $u\in S$,
then there exists no vertex of degree $1$ in $G[S]$.
\end{lem}

\pf Suppose to the contrary that there exists some $v\in S$ such
that $v$ is of degree $1$ in $G[S]$. Without loss of generality, let
$vw\in E(G)$, where $w\in S$. By Fact $2$, $v\notin D_w$. Since
$D_w\subseteq S$, $D_w$ does not dominate $v$, a contradiction. \qed

The following two lemmas, proved by Ananchuen and Plummer
\cite{PLU2}, will be used in our proof of the main theorem.

\begin{lem}  \label{lem3}
If $G$ is $3$-vertex-critical and $S$ is a cutset in $G$ such that
$\omega(G-S)\geq 4$ or $\omega(G-S)=3$, but each component has at
least $2$ vertices, then each vertex of $G-S$ is not adjacent to at
least one vertex of $S$.
\end{lem}

\begin{lem}  \label{lem4}
Let $G$ be a $3$-vertex-critical graph and suppose that $S$ is a cutset
of size $2$ in $G$, then $\omega(G-S)\leq 3$. Furthermore, if
$\omega(G-S)=3$, then $G-S$ must contain at least one singleton
component.
\end{lem}

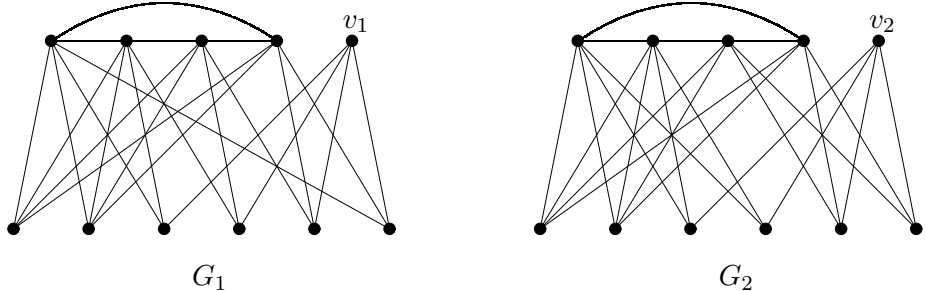
\begin{figure}[h,t]
\begin{center}
\begin{picture}(350,110)
\setlength{\unitlength}{5mm}

\put(0,1.5){\circle*{0.3}} \put(2,1.5){\circle*{0.3}}
\put(4,1.5){\circle*{0.3}} \put(6,1.5){\circle*{0.3}}
\put(8,1.5){\circle*{0.3}} \put(10,1.5){\circle*{0.3}}
\put(1,6.5){\circle*{0.3}} \put(3,6.5){\circle*{0.3}}
\put(5,6.5){\circle*{0.3}} \put(7,6.5){\circle*{0.3}}
\put(9,6.5){\circle*{0.3}} \put(4.75,0){$G_1$}
\put(8.75,6.85){$v_1$} \put(0,1.5){\line(1,5){1}}
\put(0,1.5){\line(3,5){3}} \put(0,1.5){\line(5,5){5}}
\put(0,1.5){\line(7,5){7}} \put(2,1.5){\line(-1,5){1}}
\put(2,1.5){\line(1,5){1}} \put(2,1.5){\line(3,5){3}}
\put(2,1.5){\line(5,5){5}} \put(4,1.5){\line(-3,5){3}}
\put(4,1.5){\line(-1,5){1}} \put(4,1.5){\line(5,5){5}}
\put(6,1.5){\line(-3,5){3}}\put(6,1.5){\line(-1,5){1}}
\put(6,1.5){\line(3,5){3}}
\put(8,1.5){\line(-3,5){3}}\put(8,1.5){\line(-1,5){1}}
\put(8,1.5){\line(1,5){1}}
\put(10,1.5){\line(-9,5){9}}\put(10,1.5){\line(-3,5){3}}
\put(10,1.5){\line(-1,5){1}}
\put(1,6.5){\line(2,0){2}}\put(3,6.5){\line(2,0){2}}
\put(5,6.5){\line(2,0){2}} \qbezier(1,6.5)(4,8.5)(7,6.5)

\put(14,1.5){\circle*{0.3}} \put(16,1.5){\circle*{0.3}}
\put(18,1.5){\circle*{0.3}} \put(20,1.5){\circle*{0.3}}
\put(22,1.5){\circle*{0.3}} \put(24,1.5){\circle*{0.3}}
\put(15,6.5){\circle*{0.3}} \put(17,6.5){\circle*{0.3}}
\put(19,6.5){\circle*{0.3}} \put(21,6.5){\circle*{0.3}}
\put(23,6.5){\circle*{0.3}} \put(18.75,0){$G_2$}
\put(22.75,6.85){$v_2$}
\put(14,1.5){\line(1,5){1}}\put(14,1.5){\line(3,5){3}}
\put(14,1.5){\line(5,5){5}} \put(14,1.5){\line(7,5){7}}
\put(16,1.5){\line(-1,5){1}} \put(16,1.5){\line(1,5){1}}
\put(16,1.5){\line(3,5){3}} \put(16,1.5){\line(5,5){5}}
\put(18,1.5){\line(-3,5){3}} \put(18,1.5){\line(-1,5){1}}
\put(18,1.5){\line(5,5){5}}
\put(20,1.5){\line(-5,5){5}}\put(20,1.5){\line(-3,5){3}}
\put(20,1.5){\line(3,5){3}}
\put(22,1.5){\line(-3,5){3}}\put(22,1.5){\line(-1,5){1}}
\put(22,1.5){\line(1,5){1}}
\put(24,1.5){\line(-5,5){5}}\put(24,1.5){\line(-3,5){3}}
\put(24,1.5){\line(-1,5){1}}
\put(15,6.5){\line(2,0){2}}\put(17,6.5){\line(2,0){2}}
\put(19,6.5){\line(2,0){2}} \qbezier(15,6.5)(18,8.5)(21,6.5)

\end{picture}
\caption{The graphs $G_1$ and $G_2$.} \label{fig1} \vspace{-15pt}
\end{center}
\end{figure}

Before giving our main result, we note that the graphs $G_1$ and
$G_2$ in Figure \ref{fig1} are $K_{1,5}$-free $3$-vertex-critical
$2$-connected graph of order $11$ with $\delta(G)=3$, but are not
factor-critical, since $G_i-v_i$ has no perfect matching for
$i=1,2$. We shall show that these two graphs are the only two
counterexamples for Conjecture \ref{conj}.

\begin{theo}
If $G$ is a $K_{1,5}$-free $3$-vertex-critical $2$-connected graph
of odd order with $\delta(G)\geq 3$,  except the graphs $G_1$ and
$G_2$ shown in Figure 1, then $G$ is factor-critical.
\end{theo}

\pf Suppose that $G$ is not factor-critical. By Lemma \ref{lem1} and the
parity, there exists a nonempty set $S\subseteq V(G)$ such that
$c_o(G-S)\geq |S|+1$. Without loss of generality, let $S$ be a
minimal such set with $|S|=k$. Then $k\geq 2$ as $G$ is $2$-connected.
Let $C_1,C_2,\ldots,C_t$ be the odd
components of $G-S$ and $E_1,E_2,\ldots,E_n$ the even
components of $G-S$. We consider the following cases.

\textit{Case $1$.} $k=2$.

By Lemma \ref{lem4}, then $t=3$ and $G-S$ has no even components. Since
$\delta(G)\geq 3$ and $k=2$, each odd component of $G-S$ has at
least three vertices, which contradicts to Lemma \ref{lem4}.

\textit{Case $2$.} $k=3$.

Thus, $t\geq 4$. By Lemma \ref{lem3}, each vertex of $G-S$ is not
adjacent to at least one vertex of $S$. Since $\delta(G)\geq 3$ and
$k=3$, we have $|V(C_i)|\geq 3$ for $i=1, 2, \ldots, t$. By Fact $3$, there
must exist a vertex $x$ in some odd component of $G-S$ such that
$D_x\nsubseteq S$. Clearly, $D_x\cap S\neq\emptyset$. Without loss
of generality, let $x\in V(C_1)$ and $D_x=\{u,y\}$, where $u\in S$
and $y\in V(G)-S$. Since $G$ is $K_{1,5}$-free, by the parity, so $t=4$
and $G-S$ has at most one even component.

\textbf{Claim $1$.} There exists an odd component $C_j$ ($j\geq 2$)
such that $C_j$ is a complete graph and $u$ is adjacent to every
vertex of $V(C_j)$.

If $y\in V(C_1)-\{x\}$, then $u$ is adjacent to every vertex of
$\bigcup _{i=2}^4 V(C_i)$. Since $G$ is $K_{1,5}$-free, at least two
of $C_2$, $C_3$ and $C_4$ are complete. If $y\in \bigcup _{i=2}^{4}
V(C_i)$ and suppose $y\in V(C_2)$. Then $u$ dominates all vertices
of $(V(C_1)\cup V(C_3)\cup V(C_4))-\{x\}$, and at least one of $C_3$
and $C_4$ is complete,  by $K_{1,5}$-freeness in $G$ again. If
$G-S$ has an even component $E_1$ and $y\in V(E_1)$, then $u$ is
adjacent to every vertex of $\bigcup _{i=1}^4 V(C_i)-\{x\}$. Since
$G$ is $K_{1,5}$-free, $C_2$, $C_3$ and $C_4$ are all complete. So
Claim $1$ is proved.

Without loss of generality, assume that $C_4$ is complete and $u$ is
adjacent to every vertex of $V(C_4)$.

\textbf{Claim $2$.} Each vertex of $S-\{u\}$ is not adjacent to any
vertex of $V(C_4)$.

Suppose to the contrary that $va_4\in E(G)$ for some $v\in S-\{u\}$
and $a_4\in V(C_4)$. Then $D_{a_4}\cap (\{u,v\}\cup
V(C_4))=\emptyset$, since $C_4$ is complete and $ua_4\in E(G)$. Let
$S-\{u,v\}=\{w\}$. Clearly, $w\in D_{a_4}$. Then $wa_4\notin E(G)$
and $w$ dominates $V(C_4)-\{a_4\}$. Let $b_4\in V(C_4)-\{a_4\}$.
Then $ub_4\in E(G)$ and $wb_4\in E(G)$. Consequently, $D_{b_4}\cap
(\{u,w\}\cup V(C_4))=\emptyset$ and $v\in D_{b_4}$. So $vb_4\notin
E(G)$ and $v$ dominates $V(C_4)-\{b_4\}$. Now let $c_4\in
V(C_4)-\{a_4,b_4\}$, then $c_4$ is adjacent to every vertex of $S$,
which contradicts to Lemma \ref{lem3}.

    From Claim 2, $u$ is a cut-vertex in $G$, which is against the fact that
$G$ is $2$-connected.

\textit{Case $3$.} $k=4$.

Thus, $t\geq 5$. We first show that there exists some $a\in S$ such
that $D_a\nsubseteq S$. Otherwise, $D_b\subseteq
S$ for each vertex $b\in S$. By Lemma \ref{lem2} and Fact $2$, every
vertex of $S$ in $G[S]$ has degree $0$. It is easy to check that
this is impossible.

So let $u\in S$ such that $D_u\nsubseteq S$. Clearly, $D_u\cap S\neq
\emptyset$. Let $D_u=\{v,x\}$, where $v\in S$ and $x\in V(G)-S$.
Since $G$ is $K_{1,5}$-free, so $t=5$ and $G-S$ has no even components.
Without loss of generality, let $x\in V(C_1)$, then $v$ dominates
all vertices of $\bigcup _{i=2}^{5}V(C_i)$.  Moreover,
 by $K_{1,5}$-freeness again,
$C_2$, $C_3$, $C_4$ and $C_5$ are all complete, $v$ is not adjacent
to any vertex of $V(C_1)$.

\textbf{Claim $3$.} Each vertex of $S$ is adjacent to at least three
odd components of $G-S$.

Otherwise, there exists a vertex $c\in S$ such that $c$ is
adjacent to at most two odd components of $G-S$. Let
$S'=S-\{c\}$. It is easy to see that $S'$ is a nonempty set which
satisfies the condition that $c_o(G-S')\geq |S'|+1$, contradicting
to the minimality of $S$.

Let $S-\{u,v\}=\{w,z\}$. By Claim $3$, $w$ is adjacent to at least
two of $C_2$, $C_3$, $C_4$ and $C_5$. Without loss of generality,
let $wc_i\in E(G)$, where $c_i\in V(C_i)$ for $i=2,3$. Then $z\in
D_{c_2}$. Otherwise, $u\in D_{c_2}$ and $D_{c_2}\cap V(C_1)\neq \emptyset$
 since $ux\notin E(G)$. But then $D_{c_2}$ can not
dominate $v$, a contradiction. Similarly, $z\in D_{c_3}$. Thus,
$zc_i\notin E(G)$ for $i=2,3$. By Fact $3$, then either $D_{c_2}\neq
\{u,z\}$ or $D_{c_3}\neq \{u,z\}$, say $D_{c_2}\neq \{u,z\}$. Since
$zc_3\notin E(G)$, it follows that $D_{c_2}\cap V(C_3)\neq
\emptyset$ and $z$ dominates every vertex of $V(C_1)\cup V(C_4)\cup
V(C_5)$. By similar arguments, $w\in D_{c_4}$, $w\in D_{c_5}$ for
some $c_4\in V(C_4)$ and $c_5\in V(C_5)$. Furthermore, $wc_i\notin
E(G)$ for $i=4,5$, and $w$ is adjacent to all vertices of
$V(C_1)\cup V(C_2)\cup V(C_3)$.

We next show that $C_2$ is a singleton. Otherwise, $|V(C_2)|\geq 3$ and let $a_2$, $b_2\in
V(C_2)-\{c_2\}$. By similar arguments as the above, $z\in D_{a_2}$,
$z\in D_{b_2}$ and either $D_{a_2}\neq \{u,z\}$ or $D_{b_2}\neq
\{u,z\}$. Assume that $D_{a_2}\neq \{u,z\}$. Then $D_{a_2}\cap
V(C_3)\neq \emptyset$, since $zc_3\notin E(G)$. But then $z$ is
adjacent to all vertices of $V(C_2)-\{a_2\}$ and this contradicts
to the fact that $zc_2\notin E(G)$. Similarly, $C_3$, $C_4$ and $C_5$
are all singletons of $G-S$. Since $\delta(G)\geq 3$,
$uc_i\in E(G)$ for $i=2,3,4,5$. Since $G$ is $K_{1,5}$-free, $u$ is
not adjacent to any vertex of $V(C_1)$.

Because $\delta(G)\geq 3$ and $u$, $v$ are not adjacent to any vertex
of $V(C_1)$, we have $|V(C_1)|\geq 3$. Moreover,
$D_x\cap (V(C_1)-\{x\})\neq \emptyset$ and $D_x\cap \{u,v\}\neq \emptyset$
(say, $u \in D_x$). Recall that $uv \notin E(G)$ and
$v$ is not adjacent to any vertex of $V(C_1)$, thus $v$ is not dominated by
$D_x$, a contradiction.

\textit{Case $4$.} $k=5$.

\textbf{Claim $4$.} For every vertex $x\in V(G)$, $D_x\subseteq S$.

    Otherwise, $D_u\nsubseteq S$ for some $u\in S$.
Clearly, $D_u\cap S\neq \emptyset$. Let $D_u=\{y,z\}$, where
$y\in S$ and $z\in V(G)-S$. Since $t\geq 6$, $y$ must dominate at
least $5$ odd components of $G-S$, which contradicts to the fact that
$G$ is $K_{1,5}$-free.

    Let $S=\{s_1,s_2,s_3,s_4,s_5\}$. By Fact
$3$, there are ${5\choose 2}=10$ distinct pairs of vertices in $S$
and at least $11$ vertices in $G$. So there must exist a vertex $x\in
V(G)-S$ such that $D_x\nsubseteq S$. Assume that $x\in V(C_1)$.
Clearly, $D_x\cap S\neq \emptyset$. Since $G$ is $K_{1,5}$-free,
we have $t=6$ and $G-S$ has no even components. By Claim $4$ and Lemma
\ref{lem2}, each vertex of $S$ in $G[S]$ has degree $0$ or $2$. It
is not hard to see that $G[S]$ can only be a $5$-cycle or
a union of a $4$-cycle and an isolated vertex.

\textit{Case $4.1$.} $G[S]$ is a $5$-cycle.

Let $s_1s_2s_3s_4s_5s_1$ be the $5$-cycle in the counterclockwise
order and $D_x=\{s_1,w\}$, where $w\in V(G)-S$. Since $G$ is
$K_{1,5}$-free, $w\notin V(C_1)$. Assume that $w\in V(C_2)$. Then
$s_1$ is adjacent to all vertices of $\bigcup _{i=3}^6V(C_i)$ and
$w$ dominates $s_3$, $s_4$. Moreover,  $K_{1,5}$-freeness of $G$
implies that $C_3$, $C_4$, $C_5$ and $C_6$
are all complete, $C_1$ is a singleton and $s_1$
is not adjacent to any vertex of $V(C_1)\cup V(C_2)$.

Since $D_{s_3}=\{s_1,s_5\}$, $s_5$ is adjacent to each vertex of
$V(C_1)\cup V(C_2)$. Similarly, since $D_{s_4}=\{s_1,s_2\}$, $s_2$
is adjacent to each vertex of $V(C_1)\cup V(C_2)$. Therefore, $w$ is
adjacent to all vertices of $S-\{s_1\}$. Now consider $D_w$. Since
$D_w\cap S = \{s_1\}$ and $s_1x\notin E(G)$, it follows $D_w=\{s_1,x\}$.
Hence, $x$ dominates $s_3$, $s_4$ and $V(C_2)=\{w\}$. But then
$\{s_1,s_3\}$ is a dominating set in $G$, contradicting the
assumption that $\gamma(G)=3$.

\textit{Case $4.2$.} $G[S]$ is a union of a $4$-cycle and an
isolated vertex.

Let $s_1s_2s_3s_4s_1$ be the $4$-cycle in the counterclockwise
order and $s_5$ the isolated vertex in $G[S]$. Then
$D_{s_1}=\{s_3,s_5\}$, $D_{s_2}=\{s_4,s_5\}$, $D_{s_3}=\{s_1,s_5\}$,
and $D_{s_4}=\{s_2,s_5\}$.

Since $G$ is $K_{1,5}$-free, $s_5$ is adjacent to at most $4$ odd
components of $G-S$. Without loss of generality, let
$C_1,\ldots,C_r$ be the components which are not adjacent to $s_5$.
Then $t=6$ implies $r\geq2$. Thus $s_i$ is adjacent to every vertex
of $\bigcup _{j=1}^rV(C_j)$ for $i=1,2,3,4$. Now consider $D_y$, $y
\in V(C_1)$. Clearly, $D_y \cap S = \{s_5\}$. Since $s_5$ can not
dominate $V(C_2)$, $D_y\cap V(C_2)\neq\emptyset$. Therefore, $r=2$
and $s_5$ is adjacent to every vertex of $\bigcup _{i=3}^6V(C_i)$.
Moreover, $V(C_1)=\{y\}$. By a similar argument, $C_2$ is also a
singleton.
 For each vertex $v\in\bigcup _{i=3}^6V(C_i)$, we have
$D_{v}\cap S\neq \emptyset$ and $D_{v}\nsubseteq S$, since
$s_5\notin D_{v}$ and the vertices in $S-\{s_5\}$ do not dominate
$s_5$. From $K_{1,5}$-freeness of $G$, $C_3$, $C_4$, $C_5$ and $C_6$ are
all singletons, say $V(C_i)=\{c_i\}$ for $i=3,4,5,6$.

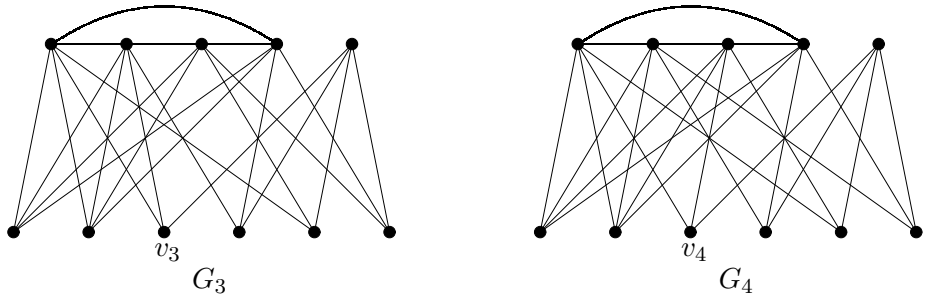
\begin{figure}[h,t]
\begin{center}
\begin{picture}(350,110)
\setlength{\unitlength}{5mm}

\put(0,1.5){\circle*{0.3}} \put(2,1.5){\circle*{0.3}}
\put(4,1.5){\circle*{0.3}} \put(6,1.5){\circle*{0.3}}
\put(8,1.5){\circle*{0.3}} \put(10,1.5){\circle*{0.3}}
\put(1,6.5){\circle*{0.3}} \put(3,6.5){\circle*{0.3}}
\put(5,6.5){\circle*{0.3}} \put(7,6.5){\circle*{0.3}}
\put(9,6.5){\circle*{0.3}} \put(4.75,0){$G_3$}
\put(3.75,0.85){$v_3$} \put(0,1.5){\line(1,5){1}}
\put(0,1.5){\line(3,5){3}} \put(0,1.5){\line(5,5){5}}
\put(0,1.5){\line(7,5){7}} \put(2,1.5){\line(-1,5){1}}
\put(2,1.5){\line(1,5){1}} \put(2,1.5){\line(3,5){3}}
\put(2,1.5){\line(5,5){5}} \put(4,1.5){\line(-3,5){3}}
\put(4,1.5){\line(-1,5){1}} \put(4,1.5){\line(5,5){5}}
\put(6,1.5){\line(-3,5){3}}\put(6,1.5){\line(1,5){1}}
\put(6,1.5){\line(3,5){3}}
\put(8,1.5){\line(-7,5){7}}\put(8,1.5){\line(-3,5){3}}
\put(8,1.5){\line(1,5){1}}
\put(10,1.5){\line(-5,5){5}}\put(10,1.5){\line(-3,5){3}}
\put(10,1.5){\line(-1,5){1}}
\put(1,6.5){\line(2,0){2}}\put(3,6.5){\line(2,0){2}}
\put(5,6.5){\line(2,0){2}} \qbezier(1,6.5)(4,8.5)(7,6.5)

\put(14,1.5){\circle*{0.3}} \put(16,1.5){\circle*{0.3}}
\put(18,1.5){\circle*{0.3}} \put(20,1.5){\circle*{0.3}}
\put(22,1.5){\circle*{0.3}} \put(24,1.5){\circle*{0.3}}
\put(15,6.5){\circle*{0.3}} \put(17,6.5){\circle*{0.3}}
\put(19,6.5){\circle*{0.3}} \put(21,6.5){\circle*{0.3}}
\put(23,6.5){\circle*{0.3}} \put(18.75,0){$G_4$}
\put(17.75,0.85){$v_4$}
\put(14,1.5){\line(1,5){1}}\put(14,1.5){\line(3,5){3}}
\put(14,1.5){\line(5,5){5}} \put(14,1.5){\line(7,5){7}}
\put(16,1.5){\line(-1,5){1}} \put(16,1.5){\line(1,5){1}}
\put(16,1.5){\line(3,5){3}} \put(16,1.5){\line(5,5){5}}
\put(18,1.5){\line(-3,5){3}} \put(18,1.5){\line(1,5){1}}
\put(18,1.5){\line(5,5){5}}
\put(20,1.5){\line(-3,5){3}}\put(20,1.5){\line(1,5){1}}
\put(20,1.5){\line(3,5){3}}
\put(22,1.5){\line(-7,5){7}}\put(22,1.5){\line(-3,5){3}}
\put(22,1.5){\line(1,5){1}}
\put(24,1.5){\line(-7,5){7}}\put(24,1.5){\line(-3,5){3}}
\put(24,1.5){\line(-1,5){1}}
\put(15,6.5){\line(2,0){2}}\put(17,6.5){\line(2,0){2}}
\put(19,6.5){\line(2,0){2}} \qbezier(15,6.5)(18,8.5)(21,6.5)

\end{picture}
\caption{The graphs $G_3$ and $G_4$.} \label{fig2} \vspace{-15pt}
\end{center}
\end{figure}


Let $H$ be the induced subgraph in $G$ with vertex set $\{s_i,c_j\
|1\leq i\leq 4,3\leq j\leq 6\}$ by deleting the edges in $G[S]$. For
$3\leq j\leq 6$, since $\delta(G)\geq 3$, $c_j$ is adjacent to at
least two vertices of $S-\{s_5\}$. On the other hand, since $G$ is
$K_{1,5}$-free, each vertex of $S-\{s_5\}$ is adjacent to at most
two vertices of $\bigcup _{i=3}^6\{c_i\}$. Thus $H$ is a $2$-regular
bipartite graph and hence consists of either a $8$-cycle or a union
of two $4$-cycles. However, there are only four such graphs under
the isomorphism, see Figure \ref{fig1} and Figure \ref{fig2}. It is
easy to see that $G_3$ and $G_4$ are not $3$-vertex-critical, since
$|D_{v_i}|> 2$ in $G_i$ for $i=3,4$. Therefore, $G_1$ and $G_2$ are
two counterexamples to Conjecture \ref{conj}.

\textit{Case $5$.} $k\geq 6$.

\textbf{Claim $5$.} For every vertex $x\in V(G)$, $D_x\subseteq S$.

Suppose that $D_x\nsubseteq S$ for some
$x\in V(G)$. Clearly, $D_x\cap S\neq \emptyset$. Let
$D_x=\{y,z\}$, where $y\in S$ and $z\in V(G)-S$. Since $t\geq 7$,
$y$ must dominate at least $5$ odd components of $G-S$,
a contradiction to $K_{1,5}$-freeness.

    Let $w$ be any vertex in $S$, then $D_w\subseteq S$ by Claim $5$.
Since $G$ is $K_{1,5}$-free, each vertex of $D_w$ can dominate at
most $4$ components of $G-S$, which implies that the number of
components of $G-S$ is at most $8$ or $t \leq 8$. That is, $6\leq
k\leq 7$.

    Let $S_i\subseteq S$ be the set of vertices in $S$ which are
adjacent to some vertex in $C_i$ for $i=1,2,\ldots,t$, and let
$d=\text{min}\{|S_i|\}$. Without loss of generality, assume that
$|S_1|=d$. Note that for any vertex $v\in V(G)-V(C_1)$, $D_v\cap
S_1\neq \emptyset$. We call such a set $D_v$ as \textit{normal
$2$-set associated with $v$ and $S_1$}, or {\it normal set} in short
. By a simple counting, we see that there are at most $k\choose
2$-$k-d\choose 2$ normal sets. Since $|V(G)-V(C_1)| \geq 2k$, Facts
3 implies $k\choose 2$-$k-d\choose 2$$\geq 2k$ or $d \geq 3$. On the
other hand, since $G$ is $K_{1,5}$-free, each vertex of $S$ is
adjacent to at most $4$ components of $G-S$, that is, $d\leq
\frac{4k}{k+1}$ or $d\leq 3$. Hence $d=3$.

\textit{Case $5.1$.} $k=6$.

Thus $t=7$ and $G-S$ has at most one even component. By Claim $5$,
there are ${6 \choose 2}=15$ distinct pairs of vertices in $S$ and
at least $13$ vertices in $G$. So by Fact $3$, $|V(G)|=13$ or $15$,
and $G-S$ has at least $6$ singletons.

It is not hard to see that there exists at least four odd components
whose corresponding $S_i$'s having the order exactly $3$, and at
least two of them are singletons. Without loss of generality, let
$C_1=\{c_1\}$ and $C_2=\{c_2\}$ be two such components. Then, for
every vertex $v\in V(G)-\{c_1\}$, $D_v\cap S_1\neq \emptyset$. There
are $12$ normal sets associated with $S_1$ in $S$, and thus
$|V(G)|=13$. Next consider $S_2$. If $S_2=S_1$, then $D_{c_2}$ can
not dominate $c_1$, a contradiction. If $|S_2\cap S_1|\leq 2$,
however, there must exist $2$ normal sets associated $S_1$ which are
not adjacent to $c_2$, at most one can be realized as $D_{c_2}$, and
the other can not dominate $c_2$, a contradiction again.

\textit{Case $5.2$.} $k=7$.

Thus $t=8$ and $G-S$ has no even components. By a similar argument
that used in the proof of Case $5.1$, one reaches the same
contradiction.

This completes the proof of our theorem.       \qed\\

\par

\noindent \title{\Large\bf Acknowledgments}
\maketitle

The authors are indebted to Professors William Y.C. Chen and Qing Cui
for bring to our attention of this conjecture and their valuable comments and discussion.

\end{document}